\documentclass [a4paper, 12pt]{article}

\usepackage{amsmath}
\usepackage{amssymb}
\usepackage{amsthm}
\usepackage{setspace}
\usepackage{multicol}

\usepackage{graphicx,color}

\renewcommand{\Re}{\mathrm{Re}}

\newtheorem{lemma}{Lemma}
\newtheorem{theorem}{Theorem}
\newtheorem*{remark}{Remark}

\numberwithin{equation}{section}
\numberwithin{lemma}{section}
\numberwithin{theorem}{section}

\begin{document}

\begin{center}
\large \textbf{Large time asymptotics for the Grinevich-Zakharov potentials}
\end{center}

\begin{center}
A. V. Kazeykina and R. G. Novikov
\end{center}

\textbf{Abstract.} In this article we show that the large time asymptotics for the Grinevich--Zakharov rational solutions of the Novikov--Veselov equation at positive energy (an analog of KdV in $ 2 + 1 $ dimensions) is given by a finite sum of localized travel waves (solitons).

\section{Introduction}
We consider the following $ 2 + 1 $--dimensional analog of the KdV equation (Novikov--Veselov equation):
\begin{equation}
\label{NV}
\begin{aligned}
& \partial_{t} v = 4 \Re \left( 4 \partial_z^3 v + \partial_z( v w ) - E \partial_z w \right), \\
& \partial_{\bar z} w = - 3 \partial_z v, \quad v = \bar v, \quad E \in \mathbb{R}, \\
& v = v( x, t ), \quad w = w( x, t ), \quad x = ( x_1, x_2 ) \in \mathbb{R}^2, \quad t \in \mathbb{R},
\end{aligned}
\end{equation}
where
\begin{equation}
\label{partials}
\partial_t = \frac{ \partial }{ \partial t }, \quad \partial_z = \frac{ 1 }{ 2 } \left( \frac{ \partial }{ \partial x_1 } - i \frac{ \partial }{ \partial x_2 } \right), \quad \partial_{ \bar z } = \frac{ 1 }{ 2 } \left( \frac{ \partial }{ \partial x_1 } + i \frac{ \partial }{ \partial x_2 } \right).
\end{equation}

We assume that
\begin{equation}
\begin{aligned}
& v \text{ is sufficiently regular and has sufficient decay as } | x | \to \infty, \\
& w \text{ is decaying as } | x | \to \infty.
\end{aligned}
\end{equation}

Equation (\ref{NV}) is contained implicitly in the paper of S.V. Manakov \cite{M} as an equation possessing the following representation:
\begin{equation}
\label{lab}
\frac{ \partial ( L - E ) }{ \partial t } = [ L - E, A ] + B( L - E )
\end{equation}
(Manakov $ L - A - B $ triple), where $ L = - \Delta + v( x, t ) $, $ \Delta = 4 \partial_z \partial_{ \bar z } $, $ A $ and $ B $ are suitable differential operators of the third and zero order respectively. Equation (\ref{NV}) was written in an explicit form by S.P. Novikov and A.P. Veselov in \cite{NV1}, \cite{NV2}, where higher analogs of (\ref{NV}) were also constructed. Note that both Kadomtsev--Petviashvili equations can be obtained from (\ref{NV}) by considering an appropriate limit $ E \to \pm \infty $ (see \cite{ZS}, \cite{G2}).

In the present article we are focused on a very interesting family of solutions for equation (\ref{NV}) for $ E = E_{ fix } > 0 $ constructed by P.G. Grinevich and V.E. Zakharov, see \cite{G1}, \cite{G2} (containing also a reference to private communication from V.E. Zakharov). The solutions of this family are given by
\begin{equation}
\label{gz_potentials}
\begin{aligned}
& v( x, t ) = - 4 \partial_z \partial_{ \bar z } \ln \det A, \\
& w( x, t ) = 12 \partial_z^2 \ln \det A, \\
\end{aligned}
\end{equation}
where $ A = (A_{lm}) $ is $ 4 N \times 4 N $--matrix,
\begin{equation}
\label{a_elements}
\begin{aligned}
& A_{ ll } = \frac{ i E^{ 1 / 2 } }{ 2 } \left( \bar z - \frac{ z }{ \lambda_l^2 } \right) - 3 i E^{ 3 / 2 } t \left( \lambda_l^2 - \frac{ 1 }{ \lambda_l^4 } \right) - \gamma_l, \\
& A_{ lm } = \frac{ 1 }{ \lambda_l - \lambda_m } \text{ for } l \neq m,
\end{aligned}
\end{equation}
$ E^{ 1 / 2 } > 0 $, $ z = x_1 + i x_2 $, $ \bar z = x_1 - i x_2 $, $ \partial_z $, $ \partial_{ \bar z } $ are defined in (\ref{partials}), and $ \lambda_1, \ldots, \lambda_{ 4 N } $, $ \gamma_1, \ldots, \gamma_{ 4N } $ are complex numbers such that
\begin{equation}
\label{lambda_conditions}
\begin{aligned}
& \lambda_j \neq 0, \quad | \lambda_j | \neq 1, \quad j = 1, \ldots, 4N, \quad \lambda_l \neq \lambda_m \text{ for } l \neq m, \\
& \lambda_{ 2 j } = - \lambda_{ 2 j - 1 }, \quad \gamma_{ 2 j - 1 } - \gamma_{ 2 j } = \frac{ 1 }{ \lambda_{ 2 j - 1 } }, \quad j = 1, \ldots, 2N, \\
& \lambda_{ 4 j - 1 } = \frac{ 1 }{ \overline{\lambda}_{ 4 j - 3 } }, \quad \gamma_{ 4 j - 1 } = \bar \lambda_{ 4 j - 3 }^2 \bar \gamma_{ 4 j - 3 }, \quad j = 1, \ldots, N.
\end{aligned}
\end{equation}

The functions $ v $, $ w $ of (\ref{gz_potentials})--(\ref{lambda_conditions}) satisfy the Novikov--Veselov equation (\ref{NV}) for positive $ E $ of (\ref{a_elements}) and have also, in particular, the following properties (see \cite{G1}, \cite{G2}):
\begin{equation}
\label{vw_properties}
\begin{aligned}
& v = \bar v, \quad w \in C^{\infty}( \mathbb{R}^2 \times \mathbb{R} ), \\
& v( x, t ), w( x, t ) \text{ are rational functions of } x \text{ and } t, \\
& v( x, t ) = O\left( | x |^{ -2 } \right), w( x, t ) = O\left( | x |^{ -2 } \right), | x | \to \infty, \text{ for each } t \in \mathbb{R};
\end{aligned}
\end{equation}
the Schr\"odinger equation $ L \psi = E \psi $, where $ L = - \Delta + v( x, t ) $,
\begin{equation}
\label{zero_scattering}
\text{ has zero scattering amplitude for fixed } E > 0 \text{ and } t \in \mathbb{R}.
\end{equation}

Because of property (\ref{zero_scattering}) the potentials $ v $ of (\ref{gz_potentials})--(\ref{lambda_conditions}) are called transparent potentials.

Note that in many respects the solutions (\ref{gz_potentials})--(\ref{lambda_conditions}) of (\ref{NV}) are similar to related solutions of the KP1 equation, see \cite{FA}, \cite{G1}, \cite{G2}.

Note that the potentials $ v $ of (\ref{gz_potentials})--(\ref{lambda_conditions}) play an important role in the direct and inverse scattering theory for the aforementioned Schr\"odinger equation $ L \psi = E \psi $ at fixed energy $ E > 0 $ and in the theory of integrable systems, see \cite{G2}, \cite{JS}, \cite{VW},  \cite{WY}, \cite{FN}, \cite{N} and references therein. However, the properties of these potentials were not yet studied sufficiently in the literature.

In order to formulate our results we use the following definition. We say that a solution $ ( v, w ) $ of (\ref{NV}) is a travel wave iff
\begin{equation}
\label{travel_wave}
v( x, t ) = V( x - c t ), \quad w( x, t ) = W( x - ct ), \quad x \in \mathbb{R}^2, \quad t \in \mathbb{R}, \\
\end{equation}
for some functions $ V $, $ W $ on $ \mathbb{R}^2 $ and some velocity $ c \in \mathbb{R}^2 $. In addition, we identify $ c = ( c_1, c_2 ) \in \mathbb{R}^2 $ with $ c = c_1 + i c_2 \in \mathbb{C} $.

The main results of the present note consist of the following:
\begin{itemize}
\item[(1)] We show that $ ( v, w ) $ of the form (\ref{gz_potentials})--(\ref{lambda_conditions}) is a travel wave iff $ N = 1 $. See Lemma 2.1 of Section \ref{main_section}.

\item[(2)] We show that there are no travel waves of the form (\ref{gz_potentials})--(\ref{lambda_conditions}), $ N = 1 $, for $ c \in \mathbb{U}_E $, and that there is an unique (modulo translations) travel wave of the form (\ref{gz_potentials})--(\ref{lambda_conditions}), $ N = 1 $, for $ c \in \mathbb{C} \backslash \mathbb{U}_E $, where $ c $ denotes travel wave velocity and $ \mathbb{U}_E $ is defined by formula (\ref{u_set}). In addition we show that there is one--to--one correspondence between permitted velocities $ c \in \mathbb{C} \backslash \mathbb{U}_E $ and $ \lambda $--sets $ \{ \lambda_1, \lambda_2, \lambda_3, \lambda_4 \} $ of (\ref{gz_potentials})--(\ref{lambda_conditions}), $ N = 1 $. See Lemma 2.2 of Section \ref{main_section}.

\item[(3)] We show that the large time asymptotics for the Grinevich--Zakharov potentials, that is for $ (v, w) $ defined by (\ref{gz_potentials})--(\ref{lambda_conditions}), is described by a sum of $ N $ localized travel waves propagating with different velocities. See Theorem \ref{main_theorem} of Section \ref{main_section}. For KP1 equation a prototype of this result was given for the first time in \cite{MZBIM}.
\end{itemize}

\section{Main results}
\label{main_section}
The main results of this article consist of Lemmas \ref{one_wave_lemma}, \ref{velocities_lemma} and Theorem \ref{main_theorem} presented below.

\begin{lemma}
\label{one_wave_lemma}
Let $ (v, w) $ be defined by (\ref{gz_potentials})--(\ref{lambda_conditions}). Then $ (v, w) $ admits the representation (\ref{travel_wave}) (and is a travel wave solution for (\ref{NV})) if and only if $ N = 1 $. In addition,
\begin{equation}
\label{c1_formula}
c = 6E \left( \bar \lambda^2 + \frac{ 1 }{ \lambda^2 } + \frac{ \lambda^2 }{ \overline{\lambda}^2 } \right)
\end{equation}
where $ c $ is the travel wave velocity and $ \lambda $ is any of $ \lambda_j $, $ j = 1, 2, 3, 4 $, which, in virtue of (\ref{lambda_conditions}), determines uniquely the $ \lambda $ set $ \{ \lambda_1, \lambda_2, \lambda_3, \lambda_4 \} $ for $ E > 0 $.
\end{lemma}
Lemma \ref{one_wave_lemma} is proved in Section \ref{proof_section}.

Let
\begin{equation}
\label{u_set}
\begin{aligned}
& \mathbb{U} = \{ u \in \mathbb{C} \colon u = r e^{ i \varphi }, r \leq | 6 ( 2 e^{ - i \varphi } + e^{ 2 i \varphi } ) |, \varphi \in [ 0, 2 \pi ] \}, \\
& \mathbb{U}_{E} = \{ u \in \mathbb{C} \colon u / E \in \mathbb{U} \}.
\end{aligned}
\end{equation}
One can see that $ \mathbb{U}_1 = \mathbb{U} $.

\begin{lemma}
\label{velocities_lemma}
\begin{itemize}
\item[(a)] Let $ c \in \mathbb{U}_E $. Then there is no travel wave solution of (\ref{NV}) of the form (\ref{gz_potentials})--(\ref{lambda_conditions}) with $ N = 1 $ and the given travel wave velocity $ c $.

\item[(b)] Let $ c \in \mathbb{C} \backslash \mathbb{U}_E $. Then there exists unique (modulo translations) solution of (\ref{NV}) of the form (\ref{gz_potentials})--(\ref{lambda_conditions}) with $ N = 1 $ and the given travel wave velocity $ c $.

\item[(c)] There is a one--to--one correspondence between $ c \in \mathbb{C} \backslash \mathbb{U}_E $ and the sets $ \{ \lambda_1, \lambda_2, \lambda_3, \lambda_4 \} $ satisfying (\ref{lambda_conditions}).
\end{itemize}
\end{lemma}
The proof of this Lemma is given in Section \ref{proof_section} and is based principally on the following auxiliary lemma.

\begin{lemma}
\label{auxiliary_lemma}
\begin{itemize}
\item[(a)] Let $ c \in \mathbb{U}_E $. Then equation (\ref{c1_formula}) has no solution $ \lambda $ satisfying $ | \lambda | \neq 1 $.
\item[(b)] Let $ c \in \mathbb{C} \backslash \mathbb{U}_E $, then equation (\ref{c1_formula}) has exactly four solutions $ \lambda_1 $, $ \lambda_2 $, $ \lambda_3 $, $ \lambda_4 $ satisfying the conditions indicated in (\ref{lambda_conditions}) for $ N = 1 $.
\end{itemize}
\end{lemma}
This Lemma is a corollary of Lemma 3.1 from \cite{KN}.

\begin{theorem}
\label{main_theorem}
Let $ ( v, w ) $ be a solution of (\ref{NV}) constructed via (\ref{gz_potentials})--(\ref{lambda_conditions}). Then the asymptotical behavior of $ ( v, w ) $ can be described as follows:
\begin{equation}
\label{vw_equivalence}
v \sim \sum_{ k = 1 }^{ N } \nu_{ k } ( \xi_k ), \quad w \sim \sum_{ k = 1 }^{ N } \omega_{ k } ( \xi_k ) \quad \text{as $ t \to \infty $},
\end{equation}
where $ \xi_k = z - c_{ 4 k } t $ and
\begin{equation}
\label{c_formula}
c_l = 6 E \left( \bar \lambda_{ l }^2 + \frac{ 1 }{ \lambda_{ l }^2 } + \frac{ \lambda_{ l }^2 }{ \overline{ \lambda }_{ l }^2 } \right).
\end{equation}
The functions $ \nu_{ k } $, $ \omega_{ k } $ are defined by the formulas
\begin{equation}
\label{nu_omega}
\begin{aligned}
& \nu_k = - 4 \partial_z \partial_{ \bar z } \ln \det A^{ ( k ) }, \\
& \omega_k = 12 \partial_z^2 \ln \det A^{ ( k ) }, \\
\end{aligned}
\end{equation}
where matrix $ A^{ ( k ) } $ is a $ 4 \times 4 $ submatrix of matrix $ A $, defined by formulas (\ref{a_elements}), such that
\begin{equation}
\label{submatrix}
A^{ ( k ) } = \{ A_{ l m } \}_{ l, m = 4 ( k - 1 ) + 1 }^{ 4 k }.
\end{equation}
\end{theorem}
\begin{remark}
The relation (\ref{vw_equivalence}) is understood in the following sense:
\begin{equation}
\label{equivalence_sense}
\lim_{ t \to \infty } v = \lim_{ t \to \infty } \sum_{ k = 1 }^{ N } \nu_k( \xi_k ) \quad \text{for fixed $ \xi = z - ct $,}
\end{equation}
where
\begin{equation}
\label{nu_sense}
\lim_{ t \to \infty } \nu_k( \xi_k ) = \begin{cases}
& 0, \quad \text{for fixed $ \xi = z - ct $, $ c \neq c_{ 4k } $}, \\
& \nu_k( \xi ), \quad \text{for fixed $ \xi = z - c_{4k} t $}.
\end{cases}
\end{equation}
\end{remark}

Theorem \ref{main_theorem} is proved in Section \ref{proof_section}. The scheme of the proof of this theorem follows principally the scheme of the derivation of the large time asymptotics for the multi--soliton solutions of the classic KdV equation (see, for example, \cite{gabov}).

\section{Proofs of Lemmas \ref{one_wave_lemma}, \ref{velocities_lemma} and Theorem \ref{main_theorem}}
\label{proof_section}
The text of the proofs presented below does not completely follow the order of statements in Section \ref{main_section} as it was constructed to form a whole logical unit. However, we specify in due course which statement is being proved.

\medskip
\noindent{\large\emph{\ref{proof_section}.1 Proof of the sufficiency part of Lemma \ref{one_wave_lemma}}}
\medskip

Let us first consider the Grinevich--Zakharov potentials defined by (\ref{gz_potentials})--(\ref{lambda_conditions}) with $ N = 1 $. Then $ A $ is a $ 4 \times 4 $ matrix and, in virtue of (\ref{lambda_conditions}), the choice of any of $ \lambda_j $, $ j = 1, 2, 3, 4 $, uniquely determines the set $ \{ \lambda_1, \lambda_2, \lambda_3, \lambda_4 \} $. Let us find $ c_j $ such that $ A_{ jj } = A_{ jj }( z - c_j t ) $. Such value $ c_j $ is a solution of the following equation
\begin{equation}
\label{equation_for_c}
\frac{ i }{ 2 } E^{ 1 / 2 } \left( \bar z - \bar c t - \frac{ 1 }{ \lambda_j^2 }( z - c t ) \right) = \frac{ i }{ 2 } E^{ 1 / 2 } \left( \bar z - \frac{ z }{ \lambda_j^2 } \right) - 3 i E^{ 3 / 2 } t \left( \lambda_j^2 - \frac{ 1 }{ \lambda_j^4 } \right).
\end{equation}
If $ | \lambda_j | \neq 1 $, then this equation is uniquely solvable and its solution is given by
\begin{equation}
\label{cj_formula}
c_j = 6 E \left( \bar \lambda^2_{ j } + \frac{ 1 }{ \lambda^2_{ j } } + \frac{ \lambda^2_{ j } }{ \bar \lambda^2_{ j } } \right).
\end{equation}
It is easy to see that due to (\ref{lambda_conditions}) $ c_1 = c_2 = c_3 = c_4 $. Thus $ A = A( z - c t ) $, and the representation (\ref{travel_wave}) with $ c $ defined by (\ref{c1_formula}) holds. Thus sufficiency in Lemma \ref{one_wave_lemma} is proved.

\medskip
\noindent{\large\emph{\ref{proof_section}.2 Proof of Lemma \ref{velocities_lemma}}}
\medskip

If $ c \in \mathbb{U}_E $, then, as follows from item (a) of Lemma \ref{auxiliary_lemma}, $ c \neq c_j $, defined by (\ref{cj_formula}), $ j = 1, 2, 3, 4 $ for any $ \lambda_1 $, $ \lambda_2 $, $ \lambda_3 $, $ \lambda_4 $ satisfying the conditions indicated in (\ref{lambda_conditions}) for $ N = 1 $. This and the sufficiency part of Lemma \ref{one_wave_lemma} imply item (a) of Lemma \ref{velocities_lemma}.

If $ c \in \mathbb{C} \backslash \mathbb{U}_E $, then, as follows from item (b) of Lemma \ref{auxiliary_lemma}, it determines via (\ref{cj_formula}) uniquely the set of $ \lambda_1 $, $ \lambda_2 $, $ \lambda_3 $, $ \lambda_4 $ satisfying the conditions indicated in (\ref{lambda_conditions}) for $ N = 1 $. Then the solution $ ( v, w ) $, constructed according to formulas (\ref{gz_potentials})--(\ref{lambda_conditions}) with $ N = 1 $, constitutes a travel wave solution of equation (\ref{NV}) with the given velocity $ c $. In the construction procedure one of the parameters $ \gamma_j $ can be chosen arbitrarily and it determines uniquely the whole set $ \{ \gamma_1, \ldots, \gamma_4 \} $.

One can see that the transform
\begin{equation*}
\begin{aligned}
& z \to z + \zeta, \\
& t \to t + \tau
\end{aligned}
\end{equation*}
turns the potential $ ( v, w ) $ into another Grinevich--Zakharov potential $ ( \tilde v, \tilde w ) $ with the parameters $ \{ \lambda_1, \ldots, \lambda_4, \tilde \gamma_1, \ldots, \tilde \gamma_4 \} $, where
\begin{equation}
\label{gammas_relation}
\gamma_j - \tilde \gamma_j = \frac{ i E }{ 2 } \left( \bar \zeta - \frac{ \zeta }{ \lambda_j^2 } \right) - 3 i E^{ 3 / 2 } \tau \left( \lambda_j^2 - \frac{ 1 }{ \lambda_j^4 } \right)
\end{equation}
for $ j = 1, 2, 3, 4 $.

On the other hand, if $ ( \tilde v, \tilde w ) $ is a Grinevich--Zakharov potential with the set of parameters $ \{ \lambda_1, \ldots, \lambda_4, \tilde \gamma_1, \ldots, \tilde \gamma_4 \} $, then it can be obtained from $ ( v, w ) $ by a translation, i.e. $ \tilde v( z, t ) = v( z + \zeta, t + \tau ) $, $ \tilde w( z, t ) = w( z + \zeta, t + \tau ) $ for appropriate $ \zeta \in \mathbb{C} $ and $ t \in \mathbb{R} $ such that (\ref{gammas_relation}) holds for some $ j $ (equations (\ref{gammas_relation}) are equivalent for $ j = 1, 2, 3, 4 $ in virtue of (\ref{lambda_conditions})). In addition, one can assume, for example, that $ \tau = 0 $ in this translation.

Thus we have proved that any $ c \in \mathbb{C} \backslash \mathbb{U}_{E} $ determines uniquely, modulo translations, the solution of (\ref{NV}) of the form (\ref{gz_potentials})--(\ref{lambda_conditions}) with $ N = 1 $ and the given travel velocity $ c $.  This proves the point (b) of Lemma \ref{velocities_lemma}.

Item (c) of Lemma \ref{velocities_lemma} follows immediately from Lemma \ref{auxiliary_lemma}. Lemma \ref{velocities_lemma} is proved.

\medskip
\noindent{\large\emph{\ref{proof_section}.3 Proof of Theorem \ref{main_theorem}}}
\medskip

Let us consider a more convenient representation of $ ( v, w ) $ defined by (\ref{gz_potentials})--(\ref{lambda_conditions}). For this purpose we first perform the differentiation with respect to $ \bar z $ in the right--hand side of formula for $ v $ in (\ref{gz_potentials}):
\begin{equation}
\label{v_repres}
v = - 4 \partial_{ z } \left[ ( \det A )^{ -1 } \partial_{ \bar z } ( \det A ) \right] = - 4 \partial_{ z } \left[ ( \det A )^{ -1 } \sum_{ i, j = 1 }^{ 4 N } \frac{ \partial A_{ i j } }{ \partial \bar z } \hat A_{ i j } \right],
\end{equation}
where $ \hat A_{ i j } $ is the $ ( i, j ) $ cofactor of the matrix $ A $.
Similarly,
\begin{equation}
\label{w_repres}
w = 12 \partial_{ z } \left[ ( \det A )^{ -1 } \sum_{ i, j = 1 }^{ 4 N } \frac{ \partial A_{ i j } }{ \partial z } \hat A_{ i j } \right].
\end{equation}

In matrix $ A $ only diagonal elements depend on $ z $, $ \bar z $, thus
\begin{equation}
\label{vw_repres}
v = - 2 i E^{ 1 / 2 } \partial_{ z } \left[ ( \det A )^{ -1 } \sum_{ j = 1 }^{ 4 N } \hat A_{ j j } \right], \quad  w = - 6 i E^{ 1 / 2 } \partial_{ z } \left[ ( \det A )^{ -1 } \sum_{ j = 1 }^{ 4 N } \frac{ 1 }{ \lambda_j^2 } \hat A_{ j j } \right].
\end{equation}

Let us consider the following families $ V^{ ( j ) } $ and $ W^{ ( j ) } $ of systems of linear algebraic equations for functions $ \psi_{ k }^{ ( j ) } $, $ \eta_{ k }^{ ( j ) } $, $ j, k = 1, \ldots 4N $:
\begin{align}
\label{v_system}
& V^{ ( j ) } \colon \quad \sum_{ k = 1 }^{ 4 N } A_{ m k } \psi_{ k }^{ ( j ) } = - 2 i E^{ 1 / 2 } \delta_{ m j } , \quad m = 1, \ldots, 4N, \\
\label{w_system}
& W^{ ( j ) } \colon \quad \sum_{ k = 1 }^{ 4 N } A_{ m k } \eta_{ k }^{ ( j ) } = - 6 i E^{ 1 / 2 } \frac{ 1 }{ \lambda_j^2 } \delta_{ m j }, \quad m = 1, \ldots, 4N.
\end{align}

Then the functions $ v $, $ w $  can be represented in the following form
\begin{equation}
\label{vw_representations}
v = \sum_{ j = 1 }^{ 4N } \frac{ \partial \psi_{ j }^{ ( j ) } }{ \partial z }, \quad w = \sum_{ j = 1 }^{ 4N } \frac{ \partial \eta_{ j }^{ ( j ) } }{ \partial z }.
\end{equation}

In order to write a system of linear algebraic equations for $ \frac{ \partial \psi_{ k }^{ ( j ) } }{ \partial z } = \left( \psi_{ k }^{ ( j ) } \right)_{ z } $, we differentiate (\ref{v_system}) with respect to $ z $:
\begin{equation*}
\sum_{ k = 1 }^{ 4 N } \left( A_{ m k } \right)_z \psi_{ k }^{ ( j ) } + \sum_{ k = 1 }^{ 4 N } A_{ m k } \left( \psi_{ k }^{ ( j ) } \right)_z = 0, \quad m = 1, \ldots, 4N,
\end{equation*}
and thus obtain
\begin{equation}
\label{v_deriv_system}
\sum_{ k = 1 }^{ 4 N } A_{ m k } \left( \psi_{ k }^{ ( j ) } \right)_z = \frac{ i E^{ 1 / 2 } }{ 2 \lambda_m^2 } \psi_{ m }^{ ( j ) }, \quad m = 1, \ldots, 4N.
\end{equation}

Now let us note that $ A_{ j j } $ can be represented in the form
\begin{equation*}
A_{ jj } = \frac{ i E^{ 1 / 2 } }{ 2 } \left[ ( \bar z - \bar c_j t ) - \frac{ 1 }{ \lambda_j^2 }( z - c_j t ) \right] - \gamma_j,
\end{equation*}
where $ c_j $ is given by formula (\ref{c_formula}). As follows from item (c) of Lemma \ref{velocities_lemma} $ c_j = c_k $ iff $ \lfloor ( j - 1 ) / 4 \rfloor = \lfloor ( k - 1 ) / 4 \rfloor $, where $ \lfloor x \rfloor $ denotes the integer part of $ x $.

Now let us fix
\begin{equation*}
\xi = z - c t
\end{equation*}
and find the limits $ \left. \psi_{ k }^{ ( j ) } \right|_{ \begin{array}{l} t \to \infty \\ \xi \text{ fixed} \end{array} } $, $ \left. \left( \psi_{ k }^{ ( j ) } \right)_{z} \right|_{ \begin{array}{l} t \to \infty \\ \xi \text{ fixed} \end{array} } $. We note that
\begin{equation*}
A_{ jj } = \frac{ i E^{ 1 / 2 } }{ 2 } \left[ \{ \bar \xi + ( \bar c - \bar c_j ) t \} - \frac{ 1 }{ \lambda_j^2 } \{ \xi + ( c - c_j ) t \} \right].
\end{equation*}
If $ c = c_j $, then $ A_{ jj } = \frac{ i E^{ 1 / 2 } }{ 2 } \left[ \bar \xi - \frac{ 1 }{ \lambda_j^2 } \xi \right] $ and is independent of $ t $. Otherwise, $ | A_{ jj } | \to \infty $ as $ t \to \infty $ at fixed $ \xi $. We substitute this into (\ref{v_system}) and consider the leading term in the Cramer's formula for $ \psi_{ k }^{ (j) } $ as $ t \to \infty $. Thus we obtain
\begin{equation*}
\begin{aligned}
& \left. \psi_{ k }^{ ( j ) } \right|_{ \begin{array}{l} t \to \infty \\ \xi \text{ fixed} \end{array} } = \hat \psi_{ k }^{ ( j ) }( \xi ), \quad k, j \colon c_k =  c_j = c, \\
& \left. \psi_{ k }^{ ( j ) } \right|_{ \begin{array}{l} t \to \infty \\ \xi \text{ fixed} \end{array} } = 0, \quad k \colon c_k \neq c \text{ or } j \colon c_j \neq c.
\end{aligned}
\end{equation*}
Here $ \hat \psi_{ k }^{ ( j ) }( \xi ) $ denotes some function of $ \xi $ independent of $ t $ at fixed $ \xi $.


Similarly, from (\ref{v_deriv_system}) we obtain that
\begin{equation*}
\begin{aligned}
& \left. \left( \psi_{ k }^{ ( j ) } \right)_{ z } \right|_{ \begin{array}{l} t \to \infty \\ \xi \text{ fixed} \end{array} } = \bar \psi_{ k }^{ ( j ) }( \xi ), \quad k, j \colon c_k =  c_j = c, \\
& \left. \left( \psi_{ k }^{ ( j ) } \right)_{ z } \right|_{ \begin{array}{l} t \to \infty \\ \xi \text{ fixed} \end{array} } = 0, \quad k \colon c_k \neq c \text{ or } j \colon c_j \neq c,
\end{aligned}
\end{equation*}
and, as previously, $ \bar \psi_{ k }^{ ( j ) }( \xi ) $ denotes some function of $ \xi $ independent of $ t $ at fixed $ \xi $.


In addition, one can see that if there exists $ k $ such that $ c = c_{ 4( k - 1 ) + 1 } = \ldots = c_{ 4k } $, then
\begin{equation}
\label{v_limit}
\left. v \right|_{ \begin{array}{l} t \to \infty \\ \xi \text{ fixed} \end{array} } = \sum_{ j = 4( k - 1 ) + 1 }^{ 4k } \bar \psi_{ j }^{ ( j ) }( \xi ) = \nu_k( \xi ),
\end{equation}
where $ \nu_k $ is defined by formula
\begin{equation}
\label{nu_k}
\nu_k = - 4 \partial_z \partial_{ \bar z } \ln \det A^{ ( k ) },
\end{equation}
matrix $ A^{ ( k ) } $ is a $ 4 \times 4 $ submatrix of matrix $ A $ from (\ref{a_elements}), such that $ A^{ ( k ) } = \{ A_{ l m } \}_{ l, m = 4( k - 1 ) + 1 }^{ 4k } $.

Similarly, for the case of function $ w $ we have
\begin{equation*}
\begin{aligned}
& \left. \eta_{ k }^{ ( j ) } \right|_{ \begin{array}{l} t \to \infty \\ \xi \text{ fixed} \end{array} } = \hat \eta_{ k }^{ ( j ) }( \xi ), \quad \left. \left( \eta_{ k }^{ ( j ) } \right)_{ z } \right|_{ \begin{array}{l} t \to \infty \\ \xi \text{ fixed} \end{array} } = \bar \eta_{ k }^{ ( j ) }( \xi ), \quad k, j \colon c_k =  c_j = c, \\
& \left. \eta_{ k }^{ ( j ) } \right|_{ \begin{array}{l} t \to \infty \\ \xi \text{ fixed} \end{array} } = 0, \quad \left. \left( \eta_{ k }^{ ( j ) } \right)_{ z } \right|_{ \begin{array}{l} t \to \infty \\ \xi \text{ fixed} \end{array} } = 0, \quad k \colon c_k \neq c \text{ or } j \colon c_j \neq c,
\end{aligned}
\end{equation*}
where $ \hat \eta_{ k }^{ j }( \xi ) $, $ \bar \eta_{ k }^{ j }( \xi ) $ are some functions of $ \xi $ independent of $ t $ at fixed $ \xi $.
%
If there exists $ k $ such that  $ c = c_{ 4( k - 1 ) + 1 } = \ldots = c_{ 4k } $, then
\begin{equation}
\label{w_limit}
\left. w \right|_{ \begin{array}{l} t \to \infty \\ \xi \text{ fixed} \end{array} } = \sum_{ j = 4( k - 1 ) + 1 }^{ 4k } \bar \eta_{ j }^{ ( j ) }( \xi ) = \omega_k( \xi ),
\end{equation}
where $ \omega_k $ is defined by formula
\begin{equation}
\label{omega_k}
\omega_k = 12 \partial_z^2 \ln \det A^{ ( k ) }
\end{equation}
and matrix $ A^{ ( k ) } $ is the same as in (\ref{nu_k}).

From (\ref{vw_representations}), (\ref{v_limit})--(\ref{omega_k}) it follows that
\begin{equation}
\label{vw_limits}
v \sim \sum_{ k = 1 }^{ N } \nu_{ k } ( \xi_k ), \quad w \sim \sum_{ k = 1 }^{ N } \omega_{ k } ( \xi_k ), \quad t \to \infty.
\end{equation}
Here $ \xi_k = z - c_{ 4 k } t $, $ \nu_k $, $ \omega_k $ are defined by (\ref{nu_omega})--(\ref{submatrix}) and the meaning of the relation (\ref{vw_limits}) is specified by (\ref{equivalence_sense})--(\ref{nu_sense}).
Theorem \ref{main_theorem} is proved.

\medskip
\noindent{\large\emph{\ref{proof_section}.4 Proof of the necessity part of Lemma \ref{one_wave_lemma}}}
\medskip

From (\ref{vw_limits}), taking into account (\ref{equivalence_sense})--(\ref{nu_sense}), one can see that $ ( v, w ) $ can be a travel wave only if $ N = 1 $. This completes the proof of Lemma \ref{one_wave_lemma}.

\noindent(A.V. Kazeykina) Centre de Math\'ematiques Appliqu\'ees, Ecole Polytechnique, Palaiseau, 91128, France \\
Lomonosov Moscow State University, GSP-1, Leninskie Gory, Moscow, 119991, Russia \\
e--mail: kazeykina@cmap.polytechnique.fr \\
(R.G. Novikov) Centre de Math\'ematiques Appliqu\'ees, Ecole Polytechnique, Palaiseau, 91128, France \\
IIEPT RAS--MITPAN, Profsoyuznaya str., 84/32, Moscow, 117997, Russia
e--mail: novikov@cmap.polytechnique.fr

\end{document}